\documentclass[12]{amsart}
\usepackage{amsmath}
\usepackage{amssymb}
\usepackage{amsfonts}
\usepackage[all, cmtip]{xy}
\newtheorem*{lemma}{Lemma}
\newtheorem*{proposition}{Proposition}
\newtheorem{definition}{Definition}
\newtheorem{corollary}{Corollary}
\newtheorem*{theorem}{Theorem}

\title  {Twisted Cyclic Homology and Crossed Product Algebras}

\author {Jack M. Shapiro}
\address {Math Department \\
Washington University \\
Saint Louis, MO  63130  USA}
\email {jshapiro@math.wustl.edu}

\begin {document}

\begin {abstract}
 $HC_*(A \rtimes G)$ is the cyclic homology of the crossed product algebra $A \rtimes G.$ For any $g \epsilon G$ we will define a homomorphism from $HC_*^g(A),$ the twisted cylic homology of $A$ with respect to $g,$ to $HC_*(A \rtimes G).$  If $G$ is the finite cyclic group generated by $g$ and $|G|=r$ is invertible in $k,$ then $HC_*(A \rtimes G)$ will be isomorphic to a direct sum of $r$ copies of $HC_*^g(A).$ For the case where $|G|$ is finite and $Q \subset k$ we will generalize the Karoubi and  Connes periodicity exact sequences for $HC_*^g(A)$  to  Karoubi and Connes periodicity exact sequences for  $HC_*(A \rtimes G)$ .
\end {abstract}

\maketitle

For an associative unitary algebra $A$ over a communtative ring $k$ together with a $k$-algebra automorphism $g$, the twisted cyclic homology of $A$ with respect to $g$, $HC_*^g(A)$, is defined in [2].  It comes from the homology of a bicomplex $(C_*, b, B)$ where $C_n$ is the quotient of the (n+1)-fold tensor product of $A$ over {k}, denoted by $A^{(n+1)}$, modulo the action of $g$. If we write $a_0 \otimes a_1\otimes \cdots \otimes a_n$ as  $(a_0 , a_1, \cdots , a_n)$, then the action of $g$ on $A^{(n+1)}$ is given by $T_g$, where $$T_g(a_0 , a_1, \cdots , a_n) = 
(g(a_0) , g(a_1), \cdots , g(a_n)).$$ Given a full (discrete) group $G$ of $k$-algebra automorphisms of $A$ we can form the crossed product algebra $A \rtimes G$  and as shown in [1, Corollary 4.2], its cyclic homology $HC_*(A \rtimes G)$ can be derived from a bicomplex $(C_*, b+ \bar b, B+T \bar B)$ with $C_n$ = $\sum_{p+q=n} k[G^{p+1}] \otimes A^{(q+1)}$. The action of $T$ on  $(g_0, \cdots, g_p / a_0, \cdots, a_q)$ $\epsilon$  $k[G^{p+1}] \otimes A^{(q+1)}$ is the identity on $k[G^{p+1}]$ and $T_g$  on $A^{(q+1)}$, where $g=g_0 \cdots g_p$. A slight alteration can be made in this bicomplex, moding out $C_n$ by the action of $T$ and dropping $T \bar B$ from the boundary maps. From this we get a homomorphism from  $HC_*^g(A)$ to $HC_*(A \rtimes G)$ for any $g \epsilon G$. For the case where $G$ is the finite cyclic group generated by $g$ with $|G|$ invertible in $k,$   $HC_*(A \rtimes G)$ will be isomorphic to a direct sum of $r$ copies of  $HC_*^g(A),$ where  $r = |G|,$ and the homomorphism will be an isomorphism of $HC_*^g(A)$ onto one of the direct summands.  It will then follow that if order($g$) = $r$, with $r$ invertible in $k,$  then  $HC_*^g(A)$ $\simeq$ $HC_*^{g^n}(A)$, for $g$ raised to any power $n$, where $(n,r)=1.$  For the case where $G$ is a finite group and $Q \subset k$ we can generalize the proceedure used in [4] for the twisted de Rham homology of $A,$ to define  $\bar HDR_*^G(A)$, the $G$-de Rham homology of $A,$ and from that get a Karoubi exact sequence $$ 0 \rightarrow \bar HDR_n^G(A) \rightarrow HC_n(A \rtimes G) \rightarrow \bar HH_{n+1}^G(A).$$ A similar thing will be done to get a Connes periodicity exact sequence. 

\section  {CYCLIC HOMOLOGIES ASSOCIATED TO $k$-ALGEBRA AUTOMORPHISMS}

A first quadrant bicomplex is a collection of $k$-modules  $C_{p,q}$  indexed by the integers  $p \ge 0$ and $q \ge 0$ together with a horizontal differential $d^h:C_{p,q} \rightarrow C_{p-1,q}$ and a vertical differential $d^v:C_{p,q} \rightarrow C_{p,q-1}$ satisfying $d^h \circ d^h=d^v \circ d^v=d^h \circ d^v+d^v \circ d^h=0$. From this we can form a chain complex, called the total complex, and its homology groups are called the homology groups of the bicomplex  [3, 1.0.11]. One such bicomplex is the one used in [3, 2.1.7.1] to define cyclic homology. It comes from maps $b:C_n \rightarrow C_{n-1}$ and $B:C_n \rightarrow C_{n+1}$ satisfying $b \circ b=B \circ B=b \circ B+B \circ b=0.$  Given a $k$-algebra automorphism, $g$, of a $k$-algebra $A$, we have the Hadfield-Kr\"ahmer bicomplex whose homology groups define the twisted cyclic homology of $A$ with respect to $g$. Here $C_n = \frac {A^{(n+1)}} {(1-T_g)}$ and the maps $b$ and $B$ are as follows.  $b = \sum_{i=0}^n (-1)^i d_i,$ where for $0 \leq i \leq {n-1},$   $d_i ( a_0, a_1,\cdots , a_n)= ( a_0,..., a_ia_{i+1},...,\cdots , a_n)$  and $d_n ( a_0, a_1,\cdots , a_n) = (g(a_n)  a_0, a_1, \cdots , a_{n-1})$.  As in cyclic homology the formula for $B$ is simpler in the normalized case where we get   $B(a_0 , a_1, \cdots , a_n)$ =  $ \sum_{i=0}^n (-1)^{ni} (1, g(a_i),\cdots , g(a_n), a_0, \cdots , a_{i-1})$ [2, section 2].

Another example of a bicomplex is the one given by Getzler-Jones in [1] to compute the cyclic homology of a crossed product algebra, $A \rtimes G$. As mentioned in the introduction $C_n = \sum_{p+q=n} k[G^{p+1}] \otimes A^{(q+1)}$, and the boundary maps are given as $b+ \bar b$ and $B + T \bar B$. The maps $b$ and $B$ on $(g_0, \cdots, g_p / a_0, \cdots, a_q)$   are the extensions of the maps $b$ and $B$ from the Hadfield-Kr\"ahmer complex with respect to $g^{-1}=(g_0 \cdots g_p)^{-1}$, that fix ${g_0, \cdots, g_p}$ in each term. The map $\bar b$ is given by $ \bar b(g_0, \cdots, g_p / a_0, \cdots, a_q) = \sum_{i=0}^{p-1}(-1)^i  (g_0, \cdots, g_i g_{i+1}, \cdots, g_p / a_0, \cdots, a_q)$ +$ (-1)^p (g_p g_0, \cdots, g_{p-1} / g_p(a_0), \cdots, g_p(a_q))$, and $\bar B$ is given by $\bar B(g_0, \cdots, g_p / a_0, \cdots, a_q)$ =
$ \sum_{i=0}^p (-1)^{ip} (1, g_{p-i+1}, \cdots, g_p, g_0, \cdots, g_{p-i} / h_i(a_0), \cdots, h_i(a_q))$, 
$h_i \equiv g_{p-i+1} \cdots g_p.$ $T(g_0, \cdots, g_p / a_0, \cdots, a_q)= (g_0, \cdots, g_p / g(a_0), \cdots, g(a_q)$, where $g$=$g_0...g_p$.  

Given a bicomplex coming from $(C_*, b, B),$ together with a $k$-module $W$, the cyclic homology with coeficients in $W$ is defined in [1, section 4] as follows. Treating $W$ as a $k[u]$-module with $deg(u)=-2$ and using the degrees in $C_*$ they form the chain complex ($C_*[[u]] \otimes_{k[u]} W, b+uB$). The homology with coeficients in $W$ is then defined as the homology of this chain complex. For a specific W we can write it in a simpler form, as pointed out in [1, sec. 4].

\begin {lemma}  For $W=k[u,u^{-1}]/u k[u]$ the homology of the chain complex  ($C_*[[u]] \otimes_{k[u]} W, b+uB$) is the  same as the homology of the bicomplex $(C_*,b,B).$ 
 \end {lemma}   

\begin {proof} For this choice of $W$ the $n^{th}$ degree term of the chain complex, $(C_*[[u]] \otimes_{k[u]} W, b+uB),$ will be $$\sum_{j=0}^{[\frac {n} {2}]} C_{n-2j} \cdot u^{-j},$$ with  boundary map $b+uB$. For the total complex of $(C_*,b,B)$ the $n^{th}$ degree term will be $$\sum_{j=0}^{[\frac {n} {2}]} C_{n-2j}$$ with boundary map is $b+B$. These two chain complexes are quasi-isomorphic.
\end {proof} 
 
From the above lemma it follows that for $W=k[u,u^{-1}]/u k[u]$ we can read corollary 4.2 in [1] as saying that the cyclic homology $HC_*(A \rtimes G)$ is isomorphic to the cyclic homolgy coming from the bicomplex   $$(\sum_{p+q=n} k[G^{p+1}] \otimes A^{(q+1)}, b+ \bar b, B + T\bar B).$$ 

\indent  We will now show that in fact for that choice of $W$ the result of corollary 4.2 can be restated as follows. 

\begin {proposition} We have an isomorphism between the cyclic homology of the crossed product algebra $HC_*(A \rtimes G)$ and the homology coming from the bicomplex $$(\sum_{p+q=n} \frac {k[G^{p+1}] \otimes A^{(q+1)}} {1- T}, b+ \bar b, B).$$  
\end {proposition}

\begin {proof} The bicomplex is well defined since it is shown in [1] that we have the equations $bB+Bb=1-T, \bar bB+B \bar b=0$,  and $T$ commutes with $b, \bar b$ and $B$. To get the $E_{pq}^1$ term in the spectral sequence used for the proof of [1, lemma 4.3], one of the main tools used is the isomorphism, $ \beta : k[G^{p+1}] \otimes A^{(q+1)} \rightarrow  k[G^p] \otimes (k[G] \otimes A^{(q+1)})$, given by $ \beta $: $(g_0, \cdots, g_p / a_0, \cdots, a_q) \rightarrow  (g_1, \cdots , g_p / g / a_0, \cdots , a_q)$, where $g=g_0 \cdots g_p$. This induces an isomorphism $$ \beta : \frac { k[G^{p+1}] \otimes A^{(q+1)}} {1- T} \rightarrow k[G^p] \otimes  \frac {k[G] \otimes A^{(q+1)}} {(1- T)},$$ which is a chain map with respect to $\bar b$ on the left,  and the boundary for group homolgy with coeficients in $ \frac {k[G] \otimes A^{(q+1)}} {(1- T)}$, on the right. Using this together with the filtration introduced in [1], after corollary 4.2, now applied to the total complex of the bicomplex $$(\sum_{p+q=n} \frac {k[G^{p+1}] \otimes A^{(q+1)}} {1- T}, b+ \bar b,  B),$$ we get that $ E_{pq}^1 = H_p(G, \frac {k[G] \otimes A^{(q+1)}} {1- \bar T})$. This is quasi-isomorphic to the $ E_{pq}^1$ term in [1] since as pointed out there, $T$ is chain homotopic to the identity on $H_p(G, k[G] \otimes A^{(q+1)})$. The rest of the proof will then be the same as in [1].
\end {proof}

\section  {RELATION BETWEEN THE HOMOLOGIES}

We can consider the bicomplex  $( \frac {k[G] \otimes A^{(n+1)}} {1- T}, b, B)$ as a sub-bicomplex of the bicomplex of the above proposition, since $\bar b$ is the zero map on $ \frac {k[G] \otimes A^{(n+1)}} {1- T}.$  The homology of this bicomplex will be denoted as $HC_*^G(A),$ and can be thought of as an extension of twisted cyclic homology to a full group, $G.$

\begin {theorem} If $G$ is a group of {k}-algebra automorphisms of a {k}-algebra $A$ then for any $g$ $\epsilon$ $G$ we have a homomorphism from the twisted cyclic homology of $A$ with respect to $g$ into $HC_*(A \rtimes G)$, $ f : HC_*^g(A) \rightarrow HC_*(A \rtimes G)$. For the case when G is the cyclic group generated by $g$, $|G| < \infty$ and $|G|$ is invertible in $k$ (e.g. ch($k$)=0), $HC_*(A \rtimes G)$ will be isomorphic to a direct sum of $r$ copies of $HC_*^g(A)$, where $r = |G|$, and f will be an isomorphism of  $HC_*^g(A)$ onto one of the summands.
\end {theorem}

\begin {proof} For any $g$ $\epsilon$ $G$, we have a map of bicomplexes, $(\frac {A^{(n+1)}} {1- T_g}, b, B) \rightarrow ( \frac {k[G] \otimes A^{(n+1)}} {1-T}, b, B)$, coming from the map sending $(a_0 , a_1, \cdots , a_n)$ to $(g^{-1}/ a_0 , a_1, \cdots , a_n),$ using the fact that $\frac {A^{(n+1)}} {1- T_g} = \frac {A^{(n+1)}} {1- T_{g^{-1}}}$.  This induces a map  $HC_*^g(A) \rightarrow HC_*^G(A).$  From the inclusion of bicomplexes $$ ( \frac {k[G] \otimes A^{(n+1)}} {1- T}, b, B) \hookrightarrow (\sum_{p+q=n} \frac {k[G^{p+1}] \otimes A^{(q+1)}} {1- T}, b + \bar b, B)$$ we get a map $HC_*^G(A) \rightarrow HC_*(A \rtimes G).$ Composing these two we get a map from $HC_*^g(A)$ to $HC_*(A \rtimes G)$, for any any $g$ $\epsilon$ $G$. If we take $G$ to be a finite group with $|G|$ invertible in $k$ then by [1, Proposition 4.6], $HC_*(A \rtimes G)$  is isomorphic to the homology of the bicomplex $( H_0(G, k[G] \otimes A^{(n+1)}), b, B) =  (\frac {k[G] \otimes A^{(n+1)}}{G},b, B)$. The action of $G$ on $ k[G] \otimes A^{(n+1)}$ involves conjugation on the $G$-factor, while $1 - T$ acts as the identity on the $G$-factor.  In that case can think of the map $HC_*^g(A) \rightarrow HC_*(A \rtimes G)$ as being induced by $$\frac {A^{(n+1)}}{1- T_g} \rightarrow  \frac {k[G] \otimes A^{(n+1)}}{1- T} \rightarrow \frac {k[G] \otimes A^{(n+1)}}{G}.$$ 

\indent  Decomposing $G$ into its cojugacy classes, $[g] = \lbrace h^{-1}gh/h \epsilon G \rbrace$, and then applying Shapiro's lemma on $H_0$, we get an isomorphism, $\frac {k[G] \otimes A^{(n+1)}}{G}$ $ \simeq$ $ \sum_{[g]} \frac {A^{(n+1)}}{G^g}$, where $G^g$ is the centralizer of $g$ in $G$ and each $A^{(n+1)}$ over $G^g$ represents the stalk of $k[G] \otimes A^{(n+1)}$ over $g$. Since for each $[g]$ the boundary maps $b$ and $B$ send the stalk to itself we have a quasi-isomorphism of bicomplexes, $(\frac {k[G] \otimes A^{(n+1)}}{G}, b, B) \simeq   (\sum_{[g]}( \frac {A^{(n+1)}}{G^g}, b, B)$. Then for each $g$ $\epsilon$ $G$ we get a composition of bicomplexes  $$( \frac {A^{(n+1)}}{1-T_g}, b, B) \rightarrow ( \frac {A^{(n+1)}}{G^g}, b, B) \hookrightarrow (\sum_{[g]} \frac {A^{(n+1)}}{G^g},b, B) \simeq (\frac {k[G] \otimes A^{(n+1)}}{G}, b, B),$$ which induces the homomorphism $HC_*^g(A) \rightarrow HC_*(A \rtimes G)$. If $G$ is the cyclic group generated by $g$, then for all $g\prime$ $\epsilon$ $G$ we have $[g\prime]$ = {$g\prime$}, $G^{g\prime}=G$ and $ \frac {A^{(n+1)}}{G}$ = $\frac {A^{(n+1)}}{g},$ which means that  $HC_*(A \rtimes G)$  is isomorphic to a direct sum of  $HC_*^g(A)$, $r$-times, where $r = |G|,$ and under the homomorphism, $HC_*^g(A)$ maps isomorphically onto one of the direct summands.
\end {proof}

\begin {corollary} If $g$ and $g \prime$ are $k$-algebra automorphisms of $A$ which generate the same finite cyclic group $G$, whose order is invertible in $k$, then $HC_*^g(A)$ $\simeq$ $HC_*^{g \prime}(A)$. In particular if $ch(k)=0$ and $|g|$ = $r$ then we have $HC_*^g(A)$ $\simeq$ $HC_*^{g^n}(A)$, for $g$ raised to any power $n, (n,r)=1$.
\end {corollary} 

 For the case where $G$ is a finite group  with $|G|$ invertible in $k$ we have the cyclic homology of $A \rtimes G$ as the homology of the bicomplex $ (\frac {k[G] \otimes A^{(n+1})}{G},b, B).$ Using the first column of that bicomplex we can extend the definition of twisted Hochschild homology given in [2] to a full group $G$.

\begin {definition}  The $G$-Hochschild homology of A, $HH_*^G(A),$ is the homology of the chain complex  $ (\frac {k[G] \otimes A^{(n+1})}{G},b),$ the first column of the bicomplex   $ (\frac {k[G] \otimes A^{(n+1})}{G},b, B).$
\end {definition}

\indent   From the proceedure used in [3, remark 2.2.2] we get from this a Connes' periodicity exact sequence, $$ ... \rightarrow HH_n^G(A) \rightarrow HC_n(A \rtimes G) \rightarrow HC_{n-2}(A \rtimes G) \rightarrow HH_{n-2}^G(A) \rightarrow ... \space .$$

\indent  Next we can generalize the proceedure used in [4] for twisted de Rham homology to define $G-$de Rham homology over a group, $G.$  First let $\bar C_n^G(A)$ $\equiv$  $\frac {k[G] \otimes A \otimes \bar A^{(n)}}{G},$ and put on it the differential defined by $d(g/ a_0, \bar a_1, ..., \bar a_n) = (g/ 1, \bar a_0, \bar a_1, ..., \bar a_n).$  Then let $\bar C_n^G(A)_{ab}$ be  $\frac {\bar C_n^G(A)}{Im(bd+db) + Im(b)}$, which is well defined since both $b$ and $d$ commute with the action of $G.$  

\begin {definition} The $G$-de Rham homology, $ \bar HDR_n^G(A),$ is the homology of the chain complex $(\bar C_*^G(A)_{ab},d).$
\end {definition}

\indent  For the case where $Q \subset k$ we automatically get that $|G|$ is invertible in $k.$ More than that we get that for a finite group $G$ the cyclic homology $HC_*(A \rtimes G)$ is isomorphic to the homology coming from the Connes complex $ ( \frac {k[G] \otimes A^{(n+1)}} {1-t}, b),$  where $t(g/a_0, ...,a_n) = (g/g^{-1}(a_n), a_0, ..., a_{n-1}).$ Using that we can  generalize the Karoubi theorem proved in [4] for twisted cyclic homology. The arguements will be the same with $\bar C_n^G(A)$  replacing  $\bar C_n^g(A).$

\begin {corollary} If $G$ is a finite group, $A$ a $k$-algebra with $Q \subset k,$  then the $G$-de Rham homology of $A$ and the cyclic homology of $A \rtimes G$ are related by the following exact sequence, $$ 0 \rightarrow \bar HDR_n^G(A) \rightarrow \bar HC_n(A \rtimes G) \rightarrow \bar HH_{n+1}^G(A).$$
\end {corollary}

\end {document}